\documentclass[12pt]{amsart}
\usepackage{amsmath}
\usepackage{amsthm}
\usepackage{amssymb}
\usepackage{amsfonts,mathrsfs}
\usepackage{stmaryrd}
\usepackage{amsxtra}
\usepackage{epsfig}
\usepackage{verbatim}
\usepackage{color}

\theoremstyle{plain}
\newtheorem{theorem}[equation]{Theorem}

\newtheorem{lemma}[equation]{Lemma}

\theoremstyle{remark}
\newtheorem{remark}{Remark}

\include{header}

\newcommand{\wB}{\mathbf{B}_{\lambda}}

\newcommand{\uD}{\mathbb{D}}
\newcommand{\wBt}{\mathbf{B}_{\lambda_t}}
\begin{document}

\bibliographystyle{alpha}

\title[Sobolev Regularity of Weighted Projections]{Sobolev Regularity of Weighted Bergman Projections on the Unit Disc}
\author{Yunus E. Zeytuncu}
\address{Department of Mathematics, Ohio State University, Columbus, Ohio 43210}
\email{yunus@math.ohio-state.edu}
\date{}
\begin{abstract} 
We show that weighted Bergman projections corresponding to radially symmetric weights on $\uD$ are bounded on Sobolev spaces.
\end{abstract}
\maketitle

\section{Introduction}
\subsection{Preliminaries}

Let $\mathbb{D}$ denote the unit disc in $\mathbb{C}^1$ and let $\lambda(r)$ be a continuous non-negative function on $[0,1)$. We consider $\lambda$ as a radial weight on $\mathbb{D}$ by setting $\lambda(z)=\lambda(|z|)$. We denote the Lebesgue measure on $\mathbb{C}$ by $dA(z)$ and  the space of square integrable functions on $\mathbb{D}$ with respect to the measure $\lambda(z)dA(z)$ by $L^2(\lambda)$. This is a Hilbert space with the inner product and the norm defined by
\begin{align*}
\left<f,g\right>_{\lambda}=\int_{\mathbb{D}}f(z)\overline{g(z)}\lambda(z)dA(z)~\text{ and }~
||f||_{\lambda}^2=\int_{\mathbb{D}}|f(z)|^2\lambda(z)dA(z).
\end{align*}

The space of holomorphic functions that are in $L^2(\lambda)$ is denoted by $A^2(\lambda)$. The Bergman inequality (see the first page of \cite{Duren04}) indicates that $A^2(\lambda)$ is a closed subspace of $L^2(\lambda)$. The orthogonal projection between these two spaces is called \textit{the weighted Bergman projection} and denoted by $\wB$, i.e.
\begin{equation*}
\mathbf{B}_{\lambda}: L^2(\lambda) \to A^2(\lambda).
\end{equation*}

It follows from the Riesz representation theorem that $\wB$ is an integral operator. The kernel is called \textit{the weighted Bergman kernel} and denoted by $B_{\lambda}(z,w)$, i.e. for any $f\in L^2(\lambda)$,
\begin{equation*}
\mathbf{B}_{\lambda}f(z)=\int_{\mathbb{D}} B_{\lambda}(z,w)f(w)\lambda(w)dA(w).
\end{equation*}

For a radial weight $\lambda$ as above, the monomials $\{z^n\}_{n=0}^{\infty}$ form an
orthogonal basis for $A^2(\lambda)$ and after normalization the weighted Bergman kernel is
given by, $B_{\lambda}(z,w)=\sum_{n=0}^{\infty}\alpha_n(z\bar w)^n$
 where $\alpha_n=\frac{1}{2\pi\int_0^1r^{2n+1}\lambda(r)dr}.$ The numbers $\alpha_n$'s are called \textit{the Bergman coefficients} of $\lambda$.

The general theory and details can be found in \cite{ForelliRudin} and \cite{Zhubook}.\\

For $k\in \mathbb{N}$, let $W^k(\lambda)$ denote the $k-$th weighted Sobolev space. The $k-$th Sobolev norm is computed as
$$||f||^2_{k,\lambda}=\sum_{|\beta|\leq k}\left|\left|\frac{\partial^{\beta}}{\partial(z,\bar{z})^{\beta}}f\right|\right|^2_{0,\lambda}=\sum_{|\beta|\leq k}\int_{\mathbb{D}}\left|\frac{\partial^{\beta}}{\partial(z,\bar{z})^{\beta}}f\right|^2\lambda(z)dA(z).$$
\vskip 1cm

\subsection{Statement}
The purpose of this paper is to show that weighted Bergman projections corresponding to radial weights are bounded on Sobolev spaces.

\begin{theorem}\label{second} Let $\lambda$ be an integrable radial weight that is non-vanishing and smooth on $\mathbb{D}$. Then $\wB$ is exactly regular i.e. $\wB$ maps $W^k(\lambda)$ to $W^k(\lambda)$ boundedly for all $k\in \mathbb{N}$. \\
\end{theorem}

\begin{remark}\label{different}
In particular, this theorem applies to weights $\lambda(r)=(1-r^2)^t$ for any $t>-1$ and $\lambda(r)=(1-r^2)^A\exp\left(\frac{-B}{(1-r^2)^{\alpha}}\right)$ for any $A\geq 0, B>0, \alpha>0$. Although, the weighted Bergman projections corresponding to these weights share the same Sobolev regularity, their $L^p$ regularity differ significantly, see \cite{Zeytuncu10}.\\
\end{remark}

\begin{remark}

The radial symmetry of the weight $\lambda$ plays a significant role in the proof. Two other places where rich symmetric structures used to prove exact regularity are \cite{Boas84} and \cite{Straube86}. The proof we present here imitates the second proof in \cite{Boas84}. \\
\end{remark}

\begin{remark}
Similar Sobolev regularity results for weighted Bergman projections appear in \cite{BonamiGrellier} and \cite{ChangLi}.\\
\end{remark}

This work is a part of my PhD dissertation at The Ohio State University. I thank my advisor, Jeffery D. McNeal, for his guidance and encouragement.\\

\section{Proof of Theorem \ref{second}}

We start with the following lemmas. 

\begin{lemma}\label{one}
For any $j\in \mathbb{N}$, there exists an operator $M_j$ and a constant $C_j>0$ such that for any two holomorphic polynomials $h$ and $g$ 
\begin{align*}\left<\frac{\partial^j}{\partial z^j}h,g\right>_{\lambda}=&\left<h,\frac{\partial^j}{\partial z^j}M_jg\right>_{\lambda}~\text{ and }~
||M_jg||_{0,\lambda}\leq C_j ||g||_{0,\lambda}.\\
\end{align*}
\end{lemma}

\begin{proof}
For any holomorphic polynomial $g(z)=\sum_{n=0}^Ng_nz^n$ and for given $j\in \mathbb{N}$, define
$$M_jg(z)=\sum_{n=0}^{N}g_n\frac{(n+j)!(n+j)!}{(n+2j)!n!}\frac{\alpha_{n+j}}{\alpha_n}z^{n+2j}$$
where $\alpha_n$'s are the Bergman coefficients of $\lambda$. We can compute the derivatives explicitly:
$$\frac{\partial^j}{\partial z^j}M_jg(z)=\sum_{n=0}^{N}g_n\frac{(n+j)!}{n!}\frac{\alpha_{n+j}}{\alpha_n}z^{n+j}.$$
We now look at the $\lambda$-inner product of this expression with a monomial $z^a$, for $a\geq j$:
\begin{align*}
\left<\frac{\partial^j}{\partial z^j}M_jg,z^a\right>_{\lambda}&=g_{a-j}\frac{a!}{(a-j)!}\frac{\alpha_a}{\alpha_{a-j}}\left<z^a,z^a\right>_{\lambda}\\
&=g_{a-j}\frac{a!}{(a-j)!}\frac{1}{\alpha_{a-j}}\\
&=g_{a-j}\frac{a!}{(a-j)!}\left<z^{a-j},z^{a-j}\right>_{\lambda}\\
&=\left<g,\frac{\partial^j}{\partial z^j}z^a\right>_{\lambda}.\\
\end{align*}
The first and the last terms are both equal to 0 for $a<j$. Hence by linearity we obtain
$$\left<g,\frac{\partial^j}{\partial z^j}h\right>_{\lambda}=\left<\frac{\partial^j}{\partial z^j}M_jg,h\right>_{\lambda}.$$
Since the weight is radial we can compute the $L^2$ norms directly from the Taylor coefficients as
\begin{align*}
||M_jg||_{0,\lambda}^2&=\sum_{n=0}^{N}|g_n|^2\left(\frac{(n+j)!(n+j)!}{(n+2j)!n!}\right)^2\frac{\alpha_{n+j}^2}{\alpha_n^2\alpha_{n+2j}}\\
||g||_{0,\lambda}^2&=\sum_{n=0}^{N}|g_n|^2\frac{1}{\alpha_n}.
\end{align*}
Noting that the sequence $\left\{\left(\frac{(n+j)!(n+j)!}{(n+2j)!n!}\right)^2\frac{\alpha_{n+j}^2}{\alpha_n\alpha_{n+2j}}\right\}_{n\in\mathbb{N}}$
is bounded (use Cauchy-Schwarz for $\alpha$ terms) for a fixed $j$, we immediately get
$$||M_jg||_{0,\lambda}\leq C_j ||g||_{0,\lambda}.$$
\end{proof}
\begin{remark}
The operator norm of $M_j$ is at most $\sup_{1\leq n< \infty}\left[\left(\frac{(n+j)!(n+j)!}{(n+2j)!n!}\right)^2\frac{\alpha_{n+j}^2}{\alpha_n\alpha_{n+2j}}\right].$ Moreover, if we keep the degree of the holomorphic polynomials $g$ less than $m$, then the operator norm of $M_j$ is at most $\sup_{1\leq n\leq m}\left[\left(\frac{(n+j)!(n+j)!}{(n+2j)!n!}\right)^2\frac{\alpha_{n+j}^2}{\alpha_n\alpha_{n+2j}}\right].$
\end{remark}

The next lemma substitutes for the holomorphic integration by parts lemma in \cite{Boas84}. Here, we assume that the weight $\lambda(r)$ vanishes at $r=1$ to infinite order. This allows us to integrate by parts without boundary terms.
\begin{lemma}\label{two}
Suppose that $\lambda(r)$ vanishes at $r=1$ to infinite order. Then for any $j\in \mathbb{N}$ there exists a constant $D_j>0$ such that for any $f\in W^j(\lambda)$  and any holomorphic polynomial $p$ we have 
\begin{equation}\label{byparts}
\left|\left<\frac{\partial^j}{\partial z^j}p,f\right>_{\lambda}\right|\leq D_j ||p||_{0,\lambda}||f||_{j,\lambda}.
\end{equation}\\
\end{lemma}

\begin{proof}
Any function $f\in W^j(\lambda)$ with support in $\{|z|<\frac{1}{2}\}$ clearly satisfies the estimate \eqref{byparts}, so we can assume that $f$ is identically zero on $\{|z|<\frac{1}{3}\}$.  

By the radial symmetry of the weight we can trade $z$ and $\bar{z}$ derivatives of $\lambda$ up to a factor. More precisely, regarding  $\lambda$ as a function of $|z|^2$ we get 
$$\frac{\partial^l}{\partial \bar{z}^l}\lambda(z)=\left(\frac{\partial^l}{\partial z^l}\lambda(z)\right)\frac{z^l}{\bar{z}^l}~\text{ for any }l\in \mathbb{N}.$$

When the support of $f$ is away from zero (and it can be even chosen away from a branch cut) we can make use of this identity. Also $\lambda$ vanishes on $b\mathbb{D}$ to infinite order, so we can integrate by parts as many times as we want without any boundary terms. Combination of these two observations with Cauchy-Schwarz inequality give the estimate
{\allowdisplaybreaks
\begin{align*}
\left|\left<\frac{\partial^j}{\partial z^j}p,f\right>_{\lambda}\right|&=\left|\left<\frac{\partial^j}{\partial z^j}p,f\lambda \right>\right|
=\left|\left<p,\frac{\partial^j}{\partial \bar{z}^j}\left(f\lambda\right) \right>\right|\\
&\leq\sum_{l+k\leq j}\left|\left<p,\frac{\partial^l}{\partial \bar{z}^l}f\frac{\partial^k}{\partial \bar{z}^k}\lambda \right>\right|\\
&=\sum_{l+k\leq j}\left|\left<p\overline{\frac{\partial^l}{\partial \bar{z}^l}f},\frac{\partial^k}{\partial \bar{z}^k}\lambda \right>\right|\\
&=\sum_{l+k\leq j}\left|\left<p\overline{\frac{\partial^l}{\partial \bar{z}^l}f},\left(\frac{\partial^k}{\partial z^k}\lambda(z)\right)\frac{z^k}{\bar{z}^k} \right>\right|\\
&=\sum_{l+k\leq j}\left|\left<p\overline{\frac{z^k}{\bar{z}^k}\frac{\partial^l}{\partial \bar{z}^l}f},\frac{\partial^k}{\partial z^k}\lambda \right>\right|\\
&=\sum_{l+k\leq j}\left|\left<\frac{\partial^k}{\partial \bar{z}^k}\left(p\overline{\frac{z^k}{\bar{z}^k}\frac{\partial^l}{\partial \bar{z}^l}f}\right),\lambda \right>\right|\\
&=\sum_{l+k\leq j}\left|\left<p\frac{\partial^k}{\partial \bar{z}^k}\left(\overline{\frac{z^k}{\bar{z}^k}\frac{\partial^l}{\partial \bar{z}^l}f}\right),\lambda \right>\right|\\
&=\sum_{l+k\leq j}\left|\left<p,\frac{\partial^k}{\partial z^k}\left(\frac{z^k}{\bar{z}^k}\frac{\partial^l}{\partial \bar{z}^l}f\right)\lambda \right>\right|\\
&=\sum_{l+k\leq j}\left|\left<p,\frac{\partial^k}{\partial z^k}\left(\frac{z^k}{\bar{z}^k}\frac{\partial^l}{\partial \bar{z}^l}f\right) \right>_{\lambda}\right|\\
&\leq D_j||p||_{0,\lambda}||f||_{j,\lambda}.
\end{align*}}
This finishes the proof of the lemma.
\end{proof}

\begin{remark}
The constant $D_j$ is independent of the weight $\lambda$.
\end{remark}

\begin{proof}[Proof of Theorem \ref{second}]
Our goal is to estimate $||\wB f||_{k,\lambda}^2=\sum_{j=0}^k||\frac{\partial^j}{\partial z^j}\wB f||^2_{0,\lambda}.$ Let $S_N$ map a holomorphic function to its $N$-th Taylor polynomial. It is clear that if we can show that for any $1\leq j\leq k$  there exists $K_j>0$ such that
\begin{equation}\label{main}
\left|\left|S_N\frac{\partial^j}{\partial z^j}\wB f\right|\right|^2_{0,\lambda}\leq K_j||f||_{k,\lambda}^2
\end{equation}
for all $N\in \mathbb{N}$ then we finish the proof.\\

{\allowdisplaybreaks
\textit{Step One.} If the weight $\lambda$ vanishes at $r=1$ to the infinite order we get this estimate directly from the lemmas above. Indeed, 
\begin{align*}
\left|\left|S_N\frac{\partial^j}{\partial z^j}\wB f\right|\right|^2_{0,\lambda}&=\sup\left\{\left|\left<h,S_N\frac{\partial^j}{\partial z^j}\wB f\right>_{\lambda}\right|~:~ \text{ for }h\in \mathcal{O}(\mathbb{D}) \text{ and }||h||_{0,\lambda}\leq1\right\}\\
&=\sup\left\{\left|\left<h,\frac{\partial^j}{\partial z^j}S_{N+j}\wB f\right>_{\lambda}\right|\dots\right\}.\\
\end{align*}}
We now concentrate on the inner product 
\begin{align*}
\left|\left<h,\frac{\partial^j}{\partial z^j}S_{N+j}\wB f\right>_{\lambda}\right|&=\left|\left<S_Nh,\frac{\partial^j}{\partial z^j}S_{N+j}\wB f\right>_{\lambda}\right| \\
&=\left|\left<\frac{\partial^j}{\partial z^j}M_jS_Nh,S_{N+j}\wB f\right>_{\lambda}\right| \text{ by the first lemma }\\
&=\left|\left<\frac{\partial^j}{\partial z^j}M_jS_Nh,\wB f\right>_{\lambda}\right|\\
&=\left|\left<\frac{\partial^j}{\partial z^j}M_jS_Nh,f\right>_{\lambda}\right|\\
&\leq D_j||M_jS_Nh||_{0,\lambda}||f||_{j,\lambda} ~\text { by the second lemma}\\
&\leq D_jC_j^N||S_Nh||_{0,\lambda}||f||_{j,\lambda} ~\text { by the first lemma}\\
&\leq D_jC_j^N||h||_{0,\lambda}||f||_{j,\lambda}. \\
\end{align*}

When we plug this estimate back into supremum calculation above we get $$\left|\left|S_N\frac{\partial^j}{\partial z^j}\wB f\right|\right|_{0,\lambda}\leq D_jC_{j,N}||f||_{j,\lambda}.$$ 
By the remarks following the lemmas, the constant $D_j$ is independent of the weight $\lambda$ and $$C_{j,N}\leq \sup_{1\leq n\leq N}\left[\left(\frac{(n+j)!(n+j)!}{(n+2j)!n!}\right)^2\frac{(\alpha_{n+j})^2}{\alpha_n\alpha_{n+2j}}\right].$$
Also we simply note that there exists $C_j>0$ such that $C_{j,N}\leq C_j$ for any $N$. This gives the desired estimate \eqref{main} for the infinite order of vanishing case.\\

\textit{Step Two.} For the weights that do not vanish to infinite order at $r=1$, we use an approximation argument. For $0<t<1$ let $\chi_t$ be a smooth radial function that is identically 1 on $\{|z|<1-t\}$, decays (without vanishing) on $\{1-t<|z|<1\}$ and vanishes to infinite order on the boundary of $\mathbb{D}$. For example, the second family of weights in Remark 1 have this property.

We set $\lambda_t=\chi_t\lambda$. Then $\lambda_t$ is still a smooth non-vanishing radial function on $\mathbb{D}$ and additionally it vanishes on $b\mathbb{D}$ at infinite order. Let $\wBt$ denote the weighted Bergman projection and $\alpha^t_n$'s denote the Bergman coefficents for the weight $\lambda_t$. By the first step we know that $\wBt$ is exactly regular.

If $f\in W^k(\lambda)$ then $f\in W^k(\lambda_t)$ for any $0<t<1$ and $||f||_{k,\lambda_t}\leq ||f||_{k,\lambda}$. A direct computation gives that for any $N,j\in \mathbb{N}$
\begin{equation}
\left|\left|S_N\frac{\partial^j}{\partial z^j}\wB f\right|\right|^2_{0,\lambda}=\lim_{t \to 0}\left|\left|S_N\frac{\partial^j}{\partial z^j}\wBt f\right|\right|^2_{0,\lambda_t}.
\end{equation}
The estimate \eqref{main} (for the case proven in the first step) implies that there exists $K_{j,N,t}>0$ such that  
$$\left|\left|S_N\frac{\partial^j}{\partial z^j}\wBt f\right|\right|^2_{0,\lambda_t}\leq K_{j,N,t}||f||_{k,\lambda_t}^2$$
where $K_{j,N,t}\leq D_j\sup_{1\leq n\leq N}\left[\left(\frac{(n+j)!(n+j)!}{(n+2j)!n!}\right)^2\frac{(\alpha_{n+j}^t)^2}{\alpha_n^t\alpha_{n+2j}^t}\right].$ For fixed $N$, if we take the limit of the previous line as $t\to 0$, then we get
$$\left|\left|S_N\frac{\partial^j}{\partial z^j}\wB f\right|\right|^2_{0,\lambda}\leq  D_j\sup_{1\leq n\leq N}\left[\left(\frac{(n+j)!(n+j)!}{(n+2j)!n!}\right)^2\frac{(\alpha_{n+j})^2}{\alpha_n\alpha_{n+2j}}\right]||f||_{k,\lambda}^2.$$
The supremum above is also finite, i.e. there exists $K_{\lambda}>0$ such that 
$$\sup_{1\leq n\leq N}\left[\left(\frac{(n+j)!(n+j)!}{(n+2j)!n!}\right)^2\frac{(\alpha_{n+j})^2}{\alpha_n\alpha_{n+2j}}\right]\leq K_{\lambda}$$
for any $N$. Hence we get
$$\left|\left|S_N\frac{\partial^j}{\partial z^j}\wB f\right|\right|^2_{0,\lambda}\leq  D_jK_{\lambda}||f||_{k,\lambda}^2$$
proving the desired estimate in the general case.\\
\end{proof}

\bibliographystyle{plain}
\bibliography{SobolevBib}
\end{document}